\documentstyle[11pt]{article}
\title{\centerline{SYMMETRY GROUP OF \c{T}I\c{T}EICA}
\centerline{SURFACES PDE}}
\author{C.UDRI\c{S}TE, N. B\^IL\u{A}}
\date{}
\begin{document}
\maketitle

\begin{abstract}

Using the theory of the symmetry groups for PDEs of order two ([7], [17], [20]),
one finds the symmetry group $G$ associated to \c{T}i\c{t}eica surfaces PDE.
One proves that Monge-Amp$\grave {\hbox{e}}$re-\c{T}i\c{t}eica
PDE which is invariant with
respect to $G'$, where $G'$ is the maximal solvable subgroup of the
symmetry group $G$, is just the PDE of \c{T}i\c{t}eica surfaces. One studies the
inverse problem and one shows that the \c{T}i\c{t}eica surfaces PDE is an
Euler-Lagrange equation. One determines the variational symmetry group
of the associated functional, and one obtains the conservation laws associated
to the \c{T}i\c{t}eica surfaces PDE. One finds some group-invariant
solutions of the \c{T}i\c{t}eica surfaces PDE.
All these results shows that \c{T}i\c{t}eica surfaces theory is strongly
related to variational problems, and hence it is a subject of global
differential geometry.

\end{abstract}



\section{Introduction}

The symmetry group (or strong symmetry group [19]) associated to a PDE
is a Lie group of local transformations which change the solutions
of the equation into its solutions. The theory of symmetry groups
has a big importance in Geometry, Mechanics and Physics ([3]-[9], [11], [13],
[14], [17]-[21], [25]).
We apply this theory in the case of \c{T}i\c{t}eica surfaces PDE in Geometry,
determining the symmetry group, some group-invariant solutions,
a \c{T}i\c{t}eica Lagrangian, and conservation laws.

The centroaffine invariant
$$
I={\frac{K}{{d^4}}}, 
$$
where $K$ is the Gauss curvature of a surface $\Sigma $ and $d$
is the distance from the origin to the tangent plane at
an arbitrary point of $\Sigma $,
was introduced by \c{T}i\c{t}eica ([22]).
A surface $\Sigma $ for which the ratio ${\frac{K}{{d^4}}}$ is constant,
is called {\it \c{T}i\c{t}eica surface}.

For the aplication of the theory of symmetry groups
([17]-[20]), we shall consider the case in which $\Sigma $
is a simple surface, being given by an explicit Cartesian
equation 
$$
\Sigma :\;\;\;u=f(x,y), 
$$
where $f\in C^2(D)$ and $D\subset {\bf R^2}$ is a domain.
In this case, the Gauss curvature of the surface $\Sigma $ is
$$
K={\frac{{u_{xx}u_{yy}-u^2_{xy}}}{{{\left(1+u^2_x+u^2_y\right)}^2}}}, 
$$
and the distance from the origin to the tangent plane at an arbitrary
point of $\Sigma $ is 
$$
d={\frac{{\vert xu_x+yu_y-u\vert }}{{\sqrt {1+u^2_x+u^2_y}}}}. 
$$
Given the nonzero function $I$ (centroaffine invariant), the condition
$$
{\frac{K}{{d^4}}}=I
$$
transcribes like a PDE 
$$
u_{xx}u_{yy}-u^2_{xy}=I{\left(xu_x+yu_y-u\right)}^4. \leqno (1) 
$$
Moreover, the restrictions $d,\;K\neq 0$ are equivalent to
$$
xu_x+yu_y-u\neq 0,\;\;u_{xx}u_{yy}-u^2_{xy}\neq 0. \leqno(2) 
$$

One proves that the symmetry group $G$ of PDE (1), with $I=constant$,
is the unimodular
subgroup of the centroaffine group.

PDE (1) is a Monge-Amp$\grave {\hbox{e}}$re equation
$$
u_{xx}u_{yy}-u^2_{xy}=H(x,y,u,u_x,u_y). \leqno(3) 
$$
Therefore PDE (3) will be called {\it Monge-Amp$\grave {\hbox{e}}$re
-\c{T}i\c{t}eica equation.}

\section{Symmetry group of a PDE of order two}

Let $D$ be an open set in {\bf R$^2$} and
$u\in C^2(D)$. The function $u^{(2)}:D\to U^{(2)}=U\times U_1\times U_2 $,
$$
u^{(2)}=(u;u_x,u_y;u_{xx},u_{xy},u_{yy}) 
$$
is called {\it the prolongation of order two of the function $u$}.

The total space $D\times U^{(2)}$ whose coordinates represent the independent
variables, the dependent variable and the derivatives of dependent variable
till the order two, is called {\it jet space of order two of the base
space $D\times U$.}

We consider the PDE of order two
$$
F(x,y,u^{(2)})=0, \leqno(4) 
$$
where $F:D\times U^{(2)}\to {\bf R}$ is a differentiable function.

{\bf Definition 1.}
PDE (4) is called of {\it maximal rank} if the associated Jacobi matrix
$$
J_F(x,y,u^{(2)})=(F_x,F_y;F_u;F_{u_x},F_{u_y};F_{u_{xx}},F_{u_{xy}},
F_{u_{yy}}) 
$$
has $rank\;1$ on the set described by the equation $F(x,y,u^{(2)})=0$.

In this case the set
$$
S=\{(x,y,u^{(2)})\in D\times U^{(2)}\vert F(x,y,u^{(2)})=0\} 
$$
is a hypersurface.

{\bf Definition 2.}
{\it The symmetry group of PDE (4)}
is a group of local transformations $G$ acting on an open set 
$M\subset D\times U$ with the properties:

(a) if $u=f(x,y)$ is a solution of the equation and if $g\cdot f$ has
sense for  $g\in G$, then $v=g\cdot f(x,y)$ is also a solution.

(b) any solution of the equation can be obtained by a DE associated
to PDE (hence any solution is $G$-invariant $g\cdot f=f,\;\forall g\in G$).

{\bf Definition 3.} Let
$$
X=\zeta (x,y,u){\frac{\partial }{{\partial x}}}+ \eta (x,y,u){\frac{\partial 
}{{\partial y}}}+\phi (x,y,u){\frac{\partial }{{\partial u}}}
$$
be a $C^{\infty }$ vector field on an open set $M\subset D\times U.$
{\it The prolongations of order one respectively two} of the vector
field $X$ are the vector fields 
$$
\hbox{pr}^{(1)}X=X+\Phi ^x{\frac{\partial }{{\partial u_x}}}+\Phi ^y {\frac{%
\partial }{{\partial u_y}}}, 
$$
$$
\hbox{pr}^{(2)}X=\hbox{pr}^{(1)}X+\Phi ^{xx}{\frac{\partial }{{\partial
u_{xx}}}}+ \Phi ^{xy}{\frac{\partial }{{\partial u_{xy}}}}+\Phi ^{yy}{\frac{%
\partial }{{\partial u_{yy}}}},\leqno (5) 
$$
where 
$$
\Phi ^x=\phi _x+(\phi _u-\zeta _{x})u_x-\eta _xu_y-\zeta_u u^2_x-\eta
_uu_xu_y, 
$$
$$
$$
$$
\Phi ^y=\phi _y-\zeta _yu_x+(\phi _u-\eta _y)u_y-\zeta _uu_xu_y
-\eta_uu^2_y
$$
and respectively
$$
\begin{array}{ccl}
\Phi^{xx} & = & \phi _{xx}+(2\phi _{xu}-\zeta _{xx})u_x-\eta _{xx}u_y+(\phi_
{uu}-2\zeta _{xu})u^2_x- \\ 
\noalign{\medskip} & - & 2\eta _{xu}u_xu_y-\zeta
_{uu}u^3_x-\eta_{uu}u^2_xu_y+ (\phi _u-2\zeta _x)u_{xx}-2\eta _xu_{xy}- \\ 
\noalign{\medskip} & - & 3\zeta _uu_xu_{xx}-\eta _uu_yu_{xx} -2\eta
_uu_xu_{xy}, 
\end{array}
$$
$$
\begin{array}{ccl}
\Phi ^{xy} & = & \phi _{xy}+(\phi _{uy}-\zeta _{xy})u_x+(\phi _{ux}-\eta
_{xy})u_y-\zeta _{uy}u^2_x+(\phi _{uu}-\zeta _{ux}- \\ 
\noalign{\medskip} & - & \eta _{uy})u_xu_y-\eta _{ux}u^2_y-\zeta
_yu_{xx}+(\phi _u-\zeta _x-\eta _y)u_{xy}-\eta _xu_{yy}- \\ 
\noalign{\medskip} & - & \zeta _uu_yu_{xx}-2\eta _uu_yu_{xy}- 2\zeta
_uu_xu_{xy}-\eta _uu_xu_{yy}- \zeta _{uu}u^2_xu_y-\eta _{uu}u_xu^2_y, 
\end{array}
$$
$$
\begin{array}{ccl}
\Phi ^{yy} & = & \phi _{yy}+(2\phi _{uy}-\eta _{yy})u_y-\zeta _{yy}u_x+
(\phi _{uu}-2\eta _{uy})u^2_y-2\zeta _{uy}u_xu_y- \\ 
\noalign{\medskip} & - & \eta _{uu}u^3_y-\zeta _{uu}u_xu^2_y+(\phi _u-2\eta
_y)u_{yy}-2\zeta _yu_{xy}- 3\eta _uu_yu_{yy}- \\ 
\noalign{\medskip} & - & \zeta _uu_xu_{yy}-2\zeta _uu_yu_{xy}. 
\end{array}
$$

For the determination of the symmetry group of PDE (4) is used
the following {\it criterion of infinitesimal invariance} [17].

{\bf Theorem 1.}
{\it Let 
$$
F(x,y,u^{(2)})=0,
$$
be a PDE of maximal rank defined on an open set $M\subset D\times U$.
If $G$ is a local group of transformations on $M$ and
$$
\hbox{pr}^{(2)}X[F(x,y,u^{(2)})]=0 \quad \hbox{whenever}\quad
F(x,y,u^{(2)})=0, \leqno (6) 
$$
for every infinitesimal generator $X$ of $G$,
then $G$ is a symmetry group of the considered equation.}

{\bf Proposition 1}.
{\it If PDE (4) defined on $M\subset D\times U$ is of maximal rank,
then the set of infinitesimal symmetries of the equation
forms a Lie algebra on $M$. Moreover, if this algebra is
finite-dimensional, then the symmetry group of PDE is a Lie
group of local transformations on $M$.}

{\bf Algorithm for determination of the symmetry group {\bf G} of PDE (4):}

-one considers the field $X$ on $M$ and its prolongations of the first
and second order, and one writes the infinitesimal invariance condition (6);

-one eliminates any dependence between partial derivatives of the function
$u$ using the given PDE;

-one writes the condition (6) like a polynomial in the partial derivatives
of $u$, and we identify this polynomial with zero;

-it follows a PDEs system in the unknown functions $\zeta,\;\eta ,\;\phi $,
and the solution of this system defines the symmetry group of PDE (4).

\section{Symmetry group of \c{T}i\c{t}eica surfaces PDE}

We consider the \c{T}i\c{t}eica surfaces PDE,
$$
u_{xx}u_{yy}-u^2_{xy}=\alpha {\left(xu_x+yu_y-u\right)}^4,\;\;\;
\alpha \in {\bf R}^{*} \leqno (1')
$$
with the conditions (2), which assure the maximal rank.

Let
$$
X=\zeta (x,y,u){\frac{\partial }{{\partial x}}}+ \eta (x,y,u){\frac{\partial 
}{{\partial y}}}+\phi (x,y,u){\frac{\partial }{{\partial u}}}
$$
be a $C^{\infty }$ vector field on the open set $M\subset D\times U$.
In the case of PDE (1$'$), the condition (6) becomes
$$
-4\alpha \zeta u_x(xu_x+yu_y-u)^3-4\alpha \eta u_y(xu_x+yu_y-u)^3+ 4\alpha
\phi (xu_x+yu_y-u)^3- 
$$
$$
-4\alpha x\Phi ^x(xu_x+yu_y-u)^3-4\alpha y \Phi ^y(xu_x+yu_y-u)^3+ \Phi
^{xx}u_{yy}- 
$$
$$
-2\Phi ^{xy}u_{xy}+\Phi ^{yy}u_{xx}=0. 
$$

Replacing the functions $\Phi ^x,\;\Phi ^y,\;\Phi ^{xx},\;\Phi ^{xy},
\;\Phi ^{yy}$ given by the relations (5) and eliminating any dependence
between partial derivatives of the function $u$ (determined
by the PDE (1$'$)), we obtain 
$$
-uu_{xx}\phi _{yy}+u_xu_{xx}(x\phi _{yy}+u\zeta _{yy})+u_yu_{xx}(y\phi
_{yy}- 2u\phi _{uy}+2u\eta _{yy})- 
$$
$$
-xu^2_xu_{xx}\zeta _{yy}+u_xu_yu_{xx}(u\zeta _{uy}-y\zeta _{yy}+2x\phi _{uy}-
x\eta _{yy})+%
$$
$$
+u^2_yu_{xx}(2y\phi _{uy}-u\phi _{uu}-y\eta _{yy}+2u\eta
_{uy})-xu^2_xu_yu_{xx} \zeta _{uy}+ 
$$
$$
+u_xu^2_yu_{xx}(u\zeta _{uu}-y\zeta _{uy}+x\phi _{uu}-2x\eta _{uy})+
u^3_yu_{xx}(y\phi _{uu}-2y\eta _{uy}+u\eta _{uu})- 
$$
$$
-xu^2_xu^2_yu_{xx}\zeta _{uu}-u_xu^3_yu_{xx}(x\eta _{uu}+y\zeta _{uu})-
yu^4_yu_{xx}\eta _{uu}+2uu_{xy}\phi _{xy}+ 
$$
$$
+2u_xu_{xy}(u\phi _{uy}-x\phi _{xy}-u\zeta _{xy})+2u_yu_{xy}(u\phi _{ux}-
y\phi _{xy}-u\eta _{xy})- 
$$
$$
-2u^2_xu_{xy}(x\phi _{uy}-x\zeta _{xy}+u\zeta _{uy})+2u_xu_yu_{xy}(u\phi
_{uu} -x\phi _{ux}-y\phi _{uy}+y\zeta _{xy}-u\zeta _{ux}+ 
$$
$$
+x\eta _{xy}-u\eta _{uy})-2u^2_yu_{xy}(y\phi _{ux}-y\eta _{xy}+u\eta _{ux})
+2xu^3_xu_{xy}\zeta _{uy}+ 2u^2_xu_yu_{xy}(y\zeta _{uy}+ 
$$
$$
+x\zeta _{ux}-u\zeta _{uu}+x\eta _{uy}-x\phi _{uu} )+2u_xu^2_yu_{xy}(y\zeta
_{ux}+y\eta _{uy}+x\eta _{ux}-u\eta _{uu}-y\phi _{uu})+%
$$
$$
+2yu^3_yu_{xy}\eta _{ux}+2xu^3_xu_yu_{xy}\zeta _{uu}+2u^2_xu^2_yu_{xy} (y\zeta
_{uu}+x\eta _{uu})+ 
$$
$$
+2yu_xu^3_yu_{xy}\eta _{uu}-uu_{yy}\phi _{xx}+u_xu_{yy}(x\phi _{xx}- 2u\phi
_{xu}+u\zeta _{xx})+ 
$$
$$
+u_yu_{yy}(x\phi _{xx}+u\eta _{xx})+u^2_xu_{yy}(2x\phi _{xu}-u\phi _{uu}
-x\zeta _{xx}+2u\zeta _{xu})+ 
$$
$$
+u_xu_yu_{yy}(2y\phi _{xu}-y\zeta _{xx}-x\eta _{xx}+2u\eta _{xu})-
yu^2_yu_{yy}\eta _{xx}+ 
$$
$$
+u^3_xu_{yy}(x\phi _{uu}-2x\zeta _{xu}+u\zeta _{uu})+u^2_xu_yu_{yy}(y\phi
_{uu} -2y\zeta _{xu}-2x\eta _{xu}+ 
$$
$$
+u\eta _{uu})-2yu_xu^2_yu_{yy}\eta _{xu}-xu^4_xu_{yy}\zeta _{uu}-
u^3_xu_yu_{yy}(y\zeta _{uu}+x\eta _{uu})- 
$$
$$
-yu^2_xu^2_yu_{yy}\eta _{uu}-2u^2_{xy}(2\phi -2x\phi _x-2y\phi _y- u(\phi
_u-\zeta _x-\eta _y))+ 
$$
$$
+u_xu^2_{xy}(4\zeta -4y\zeta _y-2x\zeta _x+2x\eta _y+2x\phi _u-4u\zeta _u)+
2u_yu^2_{xy}(2\eta -2x\eta _x- 
$$
$$
-y(\eta _y-\zeta _x-\phi _u)-2u\eta _u)+2u_{xx}u_{yy}(2\phi -2x\phi
_x-2y\phi _y-u(\phi _u-\zeta _x-\eta _y))- 
$$
$$
-2u_xu_{xx}u_{yy}(2\zeta -2y\zeta _y+x(\phi _u-\zeta _x+\eta _y)-2u\zeta
_u)- 
$$
$$
-2u_yu_{xx}u_{yy}(2\eta -2x\eta _x+y(\phi _u+\zeta _x-\eta _y)-2u\eta _u)=0. 
$$

Looking at this condition as a polynomial in the partial derivatives
of the function $u$, and identifying with the polynom zero, we obtain
the PDEs system
$$
\begin{array}{cccc}
\zeta _{xy}=0, & \zeta _{yy}=0, & \zeta _{uu}=0, & \zeta _{uy}=0, \\
\eta _{xx}=0, & \eta _{yy}=0, & \eta _{uu}=0, & \eta _{ux}=0, \\
\phi _{xx}=0, & \phi _{xy}=0, & \phi _{yy}=0, & \phi _{uy}=0,  \\
\phi _{uu}=2\eta _{uy}, & \phi _{ux}=\eta _{xy}, & \eta _{uy}=\zeta _{ux},
 & 2\phi _{xu}=\zeta _{xx}, 
\end{array}
$$
$$
2\phi -2x\phi _x-2y\phi _y-u(\phi _u-\zeta _x-\eta _y)=0,
$$
$$
2\zeta -2y\zeta _y-x(\zeta _x-\eta _y-\phi _u)-2u\zeta _u=0, 
$$
$$
2\eta -2x\eta _x+y(\phi _u+\zeta _x-\eta _y)-2u\eta _u=0, 
$$
whose solution defines the symmetry group of the equation (1$'$).
The general solution of this PDEs system is
$$
\left\{
\begin{array}{ccl}
\zeta (x,y,u) & = & C_1x+C_3y+C_4u,\\
\eta (x,y,u) & = & C_5x+C_2y+C_6u, \\ 
\phi (x,y,u) & = & C_7x+C_8y-(C_1+C_2)u,
\end{array} \right.
\leqno(7)
$$
where $C_1,...,C_8\in ${\bf R}, and consequently the infinitesimal generator of the
symmetry group $G$ is
$$
X=C_1\left(x{\frac{\partial }{{\partial x}}}-u{\frac{\partial }{{\partial u}}%
} \right)+C_2\left(y{\frac{\partial }{{\partial y}}}-u {\frac{\partial }{{%
\partial u}}}\right)+C_3y{\frac{\partial }{{\partial x}}}+ C_4u{\frac{%
\partial }{{\partial x}}}+ 
$$
$$
+C_5x{\frac{\partial }{{\partial y}}}+C_6u{\frac{\partial }{{\partial y}}}+
C_7x{\frac{\partial }{{\partial u}}}+C_8y{\frac{\partial }{{\partial u}}}. 
$$

{\bf Theorem 2.}
{\it The Lie algebra {\bf g} of the symmetry group $G$
associated to \c{T}i\c{t}eica surfaces PDE is generated by the vector fields
$$
X_1=x{\frac{\partial }{{\partial x}}}-u{\frac{\partial }{{\partial u}}},\;\;
X_2=y{\frac{\partial }{{\partial y}}}-u{\frac{\partial }{{\partial u}}},\;\;
X_3=y{\frac{\partial }{{\partial x}}},\;\;X_4=u{\frac{\partial }{{\partial x}%
}}\leqno (8)
$$
$$
X_5=x{\frac{\partial }{{\partial y}}},\;\;X_6=u{\frac{\partial }{{\partial y}%
}},\;\; X_7=x{\frac{\partial }{{\partial u}}},\;\;X_8=y{\frac{\partial }{{%
\partial u}}}, %
$$
and $G$ is the unimodular subgroup of centroaffine group.}

The constants of the structure of the
Lie algebra of the group $G$ are finding from
the table

\vspace{0,5cm}
\centerline{\begin{tabular}{|c|c|c|c|c|c|c|c|c|}
\hline
$\scriptstyle{[.,.]}$ & $\scriptstyle{X_1}$ & $\scriptstyle{X_2}$ &
$\scriptstyle{X_3}$ & $\scriptstyle{X_4}$ & $\scriptstyle{X_5}$ &
$\scriptstyle{X_6}$ & $\scriptstyle{X_7}$ & $\scriptstyle{X_8}$ \\
\hline
$\scriptstyle{X_1}$ & $\scriptstyle{0}$ & $\scriptstyle{0}$ &
$\scriptstyle{-X_3}$ & $\scriptstyle{-2X_4}$ & $\scriptstyle{X_5}$ &
$\scriptstyle{-X_6}$ & $\scriptstyle{2X_7}$ & $\scriptstyle{X_8}$ \\
\hline
$\scriptstyle{X_2}$ & $\scriptstyle{0}$ & $\scriptstyle{0}$ &
$\scriptstyle{X_3}$ & $\scriptstyle{-X_4}$ & $\scriptstyle{-X_5}$ &
$\scriptstyle{-2X_6}$ & $\scriptstyle{X_7}$ & $\scriptstyle{2X_8}$ \\
\hline
$\scriptstyle{X_3}$ & $\scriptstyle{X_3}$ & $\scriptstyle{-X_3}$ &
$\scriptstyle{0}$ & $\scriptstyle{0}$ &
$\scriptstyle{X_2-X_1}$ & $\scriptstyle{-X_4}$ & $\scriptstyle{X_8}$ &
$\scriptstyle{0}$ \\
\hline
$\scriptstyle{X_4}$ & $\scriptstyle{2X_4}$ & $\scriptstyle{X_4}$ &
$\scriptstyle{0}$ & $\scriptstyle{0}$ & $\scriptstyle{X_6}$ &
$\scriptstyle{0}$ & $\scriptstyle{-X_1}$ & $\scriptstyle{-X_3}$  \\
\hline
$\scriptstyle{X_5}$ & $\scriptstyle{-X_5}$ & $\scriptstyle{X_5}$ &
$\scriptstyle{X_1-X_2}$ & $\scriptstyle{-X_6}$ & $\scriptstyle{0}$ &
$\scriptstyle{0}$ & $\scriptstyle{0}$ & $\scriptstyle{X_7}$ \\
\hline
$\scriptstyle{X_6}$ & $\scriptstyle{X_6}$ & $\scriptstyle{2X_6}$ &
$\scriptstyle{X_4}$ & $\scriptstyle{0}$ & $\scriptstyle{0}$ &
$\scriptstyle{0}$ & $\scriptstyle{-X_5}$ & $\scriptstyle{-X_2}$ \\
\hline
$\scriptstyle{X_7}$ & $\scriptstyle{-2X_7}$ & $\scriptstyle{-X_7}$ &
$\scriptstyle{-X_8}$ & $\scriptstyle{X_1}$ & $\scriptstyle{0}$ &
$\scriptstyle{X_5}$ & $\scriptstyle{0}$ & $\scriptstyle{0}$ \\
\hline
$\scriptstyle{X_8}$ & $\scriptstyle{-X_8}$ & $\scriptstyle{-2X_8}$ &
$\scriptstyle{0}$ & $\scriptstyle{X_3}$ & $\scriptstyle{-X_7}$ &
$\scriptstyle{X_2}$ & $\scriptstyle{0}$ & $\scriptstyle{0}$ \\
\hline
\end{tabular}}
\vspace{0,5cm}                                                          

Now we shall study the converse of the Theorem 2: given the Lie group $G$
of transformations, determine the most general
Monge-Amp$\grave {\hbox{e}}$re-\c{T}i\c{t}eica PDE which admits
$G$ like group of symmetries. This implies the using of a maximal chain
of Lie subalgebras of the algebra {\bf g} of the group $G$, in the case
in which {\bf g} is solvable.

Since the Lie algebra {\bf g} of the symmetry group $G$ is not solvable,
one considers the maximal solvable Lie subalgebra {\bf g$^{\prime}$},
described by the vector fields $X_1,X_2,X_3,X_7$.
Denote $G^{\prime}\subset G$ the corresponding subgroup.

{\bf Theorem 3.}
{\it The PDE of type Monge-Amp$\grave {\hbox{e}}$re-\c{T}i\c{t}eica
of maximal rank, which admits $G^{\prime}$ like group of symmetry, is a PDE
of type \c{T}i\c{t}eica.}

{\bf Proof.}
We consider the maximal chain of Lie subalgebras of the Lie algebra
{\bf g$^{\prime}$},
$$
\{X_8\}\subset \{X_3,X_8\}\subset \{X_3,X_7\}\subset \{X_1,X_3,X_7\}\subset
\{X_1,X_2,X_3,X_7\}.
\leqno(9)
$$
We impose the condition that PDE (3) to be invariant with respect
to every of these subalgebras, denoting
$$
F=u_{xx}u_{yy}-u^2_{xy}-H(x,y,u,u_x,u_y). 
$$

\par
\noindent
1) We start with $\{X_8\}$: $X_8=y{\frac{\partial }{{\partial u}}}$
and $pr^{(2)}X_8= y{\frac{\partial }{{\partial u}}}+{\frac{\partial }
{{\partial u_y}}}.$
\par
\noindent
The condition (6) implies $pr^{(2)}X_8(F)=0$.
It follows 
$$
F=u_{xx}u_{yy}-u^2_{xy}-H_1(x,y,u_x,yu_y-u). 
$$
\par
\noindent
2) If we use $\{X_3,X_8\}$: $X_3=y{\frac{\partial }{{\partial x}}}$ and
$$
pr^{(2)}X_3 =y{\frac{\partial }{{\partial x}}}-u_x{\frac{\partial }{{%
\partial u_y}}}- u_{xx}{\frac{\partial }{{\partial u_{xy}}}} -2u_{xy}{\frac{%
\partial }{{\partial u_{yy}}}}, 
$$
then we obtain 
$$
F=u_{xx}u_{yy}-u^2_{xy}-H_2(y,u,u_x,xu_x+yu_y-u). 
$$
\par
\noindent
3) For $\{X_3,X_7\}$: $X_7=x{\frac{\partial }{{\partial u}}}$ and $pr^{(2)}X_7
=x{\frac{\partial }{{\partial u}}}+{\frac{\partial }{{\partial u_x}}},$
\par
\noindent
we find 
$$
F=u_{xx}u_{yy}-u^2_{xy}-H_3(y,xu_x+yu_y-u). 
$$
\par
\noindent
4) For $\{X_1,X_3,X_7\}$: $X_1=x{\frac{\partial }{{\partial x}}}-u{\frac{%
\partial }{{\partial u}}}$, with
$$
pr^{(2)}X_1=x{\frac{\partial }{{\partial x}}}-u{\frac{\partial }{{\partial u}%
}} -2u_x{\frac{\partial }{{\partial u_x}}}-u_y{\frac{\partial }{{\partial u_y%
}}}- 3u_{xx}{\frac{\partial }{{\partial u_{xx}}}}-2u_{xy}{\frac{\partial }{{%
\partial u_{xy}}}}-u_{yy}{\frac{\partial }{{\partial u_{yy}}}},
$$
we get 
$$
F=u_{xx}u_{yy}-u^2_{xy}-(xu_x+yu_y-u)^4H_4(y). 
$$
\par
\noindent
5) Finally, $\{X_1,X_2,X_3,X_7\}$: $X_2=y{\frac{\partial }{{\partial y}}}-
u{\frac{\partial }{{\partial u}}}$ and 
$$
pr^{(2)}X_2=y{\frac{\partial }{{\partial y}}}-u{\frac{\partial }{{\partial u}%
}}- u_x{\frac{\partial }{{\partial u_x}}}-2u_y{\frac{\partial }{{\partial u_y%
}}}- u_{xx}{\frac{\partial }{{\partial u_{xx}}}} -2u_{xy}{\frac{\partial }{{%
\partial u_{xy}}}}-3u_{yy}{\frac{\partial }{{\partial u_{yy}}}},
$$
imply 
$$
F=u_{xx}u_{yy}-u^2_{xy}-\alpha (xu_x+yu_y-u)^4,\;\;\alpha \in {\bf R}, 
$$
and consequently the Monge-Amp$\grave {\hbox{e}}$re-\c{T}i\c{t}eica PDE
is reduced to \c{T}i\c{t}eica PDE
$$
u_{xx}u_{yy}-u^2_{xy}=\alpha (xu_x+yu_y-u)^4,\;\;\alpha \in {\bf R}. 
$$
If $\alpha \neq 0$, then the condition of maximal rank is satisfied.

\section{Inverse problem associated to a PDE}

The simple form of the inverse problem in the calculus of variations
is to determine if an operator with partial derivatives is identically to
an Euler-Langrange operator with partial derivatives ([1], [2], [12], [17],
[20]). We quote

{\bf Theorem 4.}
{\it Let $T$ be the operator associated to PDE (4). $T$ is identically to
an Euler-Lagrange operator if and only if the integrability Helmholtz
conditions}
$$
\left\{
\begin{array}{ccl}
{\frac{\partial T}{{\partial u_x}}} & = & D_x\left(
{\frac{\partial T}{{\partial u_{xx}}}}\right)+D_y\left({\frac{1}{{2}}}{\frac{%
\partial T}{{\partial u_{xy}}}} \right) \\ {\frac{\partial T}{{\partial u_y}}%
} & = & D_x\left(
{\frac{1}{{2}}}{\frac{\partial T}{{\partial u_{xy}}}}\right)+D_y\left({\frac{%
\partial T}{{\partial u_{yy}}}} \right),  
\end{array}
\right.\leqno(10) 
$$
{\it are satisfied, where} 
$$
\begin{array}{ccl}
D_x & = & {\frac{\partial }{{\partial x}}}+u_x{\frac{\partial }{{\partial u}}%
}+u_{xx} {\frac{\partial }{{\partial u_x}}}+u_{xy}{\frac{\partial }{{%
\partial u_y}}}, \\ D_y & = & {\frac{\partial }{{\partial y}}}+u_y{\frac{%
\partial }{{\partial u}}}+u_{xy} {\frac{\partial }{{\partial u_x}}}+u_{yy}{
\frac{\partial }{{\partial u_y}}}. 
\end{array}
\leqno(11) 
$$

In this case, there exists a Lagrangian $L$ such that the Euler-Lagrange PDE
$E(L)(u)=0$ is equivalent to the PDE associated to the operator $T$, in the
sense that every solution of the equation $T(u)=0$ is a solution of the
Euler-Lagrange equation $E(L)(u)=0$ and conversely.

For the associated Lagrangian of order two
$$
L=L(x,y,u,u_x,u_y,u_{xx},u_{xy},u_{yy}), 
$$
{\it the Euler-Lagrange operator of order two} is 
$$
\begin{array}{ccl}
E(L)(u) & = & {\frac{\partial L}{{\partial u}}}-D_x\left({\frac{\partial L}{{%
\partial u_x}}} \right)-D_y\left({\frac{\partial L}{{\partial u_y}}}\right)+
\\ \noalign{\medskip} & + & D_{xx}\left({\frac{\partial L}{{\partial u_{xx}}}%
}\right)+D_{xy}\left( {\frac{\partial L}{{\partial u_{xy}}}}%
\right)+D_{yy}\left({\frac{\partial L}{{\partial u_{yy}}}}\right). 
\end{array}
\leqno(12) 
$$

{\bf Definition 4}.
{\it An operator $T$ is equivalent to an Euler-Lagrange operator $E(L)$,}
if there exists a nonzero function $f=f(x,y,u,u_x,u_y)$ such that
$f\cdot T=E(L)$. The function $f$ is called
{\it variational integrant factor}.

\section{Lagrangians associated to \c{T}i\c{t}eica surfaces PDE}

We consider the PDE of type \c{T}i\c{t}eica (1$'$) under conditions (2).
The operator
$$
T(u)=u_{xx}u_{yy}-u^2_{xy}-\alpha (xu_x+yu_y-u)^4,\;\;\alpha \in {\bf R}^{*},
$$
which defines the equation (1$'$), is not identically to an 
Euler-Lagrange operator, since the integrability conditions (10) are not
satisfied.

{\bf Theorem 5.}
{\it The operator $T$ is equivalent to an Euler-Lagrange operator.}

{\bf Proof.}
Suppose there exists a variational integrant factor,
$$
f=f(x,y,u,u_x,u_y),
$$
such that $f\cdot T=E(L)$.
In this case, the integrability conditions (10), for $f\cdot T$,
become 
$$
\left\{
\begin{array}{ccl}
u_{yy}\left({\frac{\partial f}{{\partial x}}}+u_x{\frac{\partial f}{{%
\partial u}}} \right)-u_{xy}\left({\frac{\partial f}{{\partial y}}}+u_y{
\frac{\partial f}{{\partial u}}}\right)+ &  &  \\  
&  &  \\ 
+\alpha {\frac{\partial f}{{\partial u_x}}}\left(xu_x+yu_y-u\right)^4+
4\alpha xf\left(xu_x+yu_y-u\right)^3 & = & 0 \\  
&  &  \\ 
u_{xx}\left({\frac{\partial f}{{\partial y}}}+u_y{\frac{\partial f}{{%
\partial u}}} \right)-u_{xy}\left({\frac{\partial f}{{\partial x}}}+u_x{
\frac{\partial f}{{\partial u}}}\right)+ &  &  \\  
&  &  \\ 
+\alpha {\frac{\partial f}{{\partial u_y}}}\left(xu_x+yu_y-u\right)^4+
4\alpha yf\left(xu_x+yu_y-u\right)^3 & = & 0. 
\end{array}
\right. 
$$

Equating to zero the coefficients of partial derivatives of
second order of the function $u$, we obtain the following PDEs system
$$
\left\{
\begin{array}{ccl}
{\frac{\partial f}{{\partial x}}}+u_x{\frac{\partial f}{{\partial u}}} & = & 
0 \\  
&  &  \\ 
{\frac{\partial f}{{\partial y}}}+u_y{\frac{\partial f}{{\partial u}}} & = & 
0 \\  
&  &  \\ 
(xu_x+yu_y-u){\frac{\partial f}{{\partial u_x}}}+4xf & = & 0 \\  
&  &  \\ 
(xu_x+yu_y-u){\frac{\partial f}{{\partial u_y}}}+4yf & = & 0.
\end{array}
\right. 
$$
The solution of this system is
$$
f(x,y,u,u_x,u_y)={\frac{C}{{(xu_x+yu_y-u)^4}}},\;\;C\in {\bf R^{*}}. 
$$
Hence, PDE of \c{T}i\c{t}eica surfaces, written in the initial form
$$
{\frac{K}{{d^4}}}=\alpha, 
$$
is an Euler-Lagrange equation.

{\bf Theorem 6.}
{\it A Lagrangian of order two associated to \c{T}i\c{t}eica
surfaces PDE is}
$$
L(x,y,u^{(2)})={\frac{u(u^2_{xy}-u_{xx}u_{yy})}{{(xu_x+yu_y-u)^4}}}- \alpha
u.\leqno(13) 
$$

{\bf Proof.}
Using formula (12), after tedious computations it follows
$$
E(L)(u)={\frac{u_{xx}u_{yy}-u^2_{xy}}{{(xu_x+yu_y-u)^4}}}-\alpha . 
$$

\section{Variational symmetry group. Conservation laws}

We will make a short presentation of the variational symmetry group
([17], [20]) for the functionals of the form
$$
{\cal L}[u]=\int _{}^{}\int _{D _0}^{}L(x,y,u^{(2)})dxdy,
$$
where $D _0$ is a domain in {\bf R$^2$}.

Let $D\subset D_0$ be a subdomain, $U$ an open set in {\bf R}
and $M\subset D\times U$ an open set. Let $u\in C^2(D),\;u=f(x,y)$
such that
$$
\Gamma _u=\{(x,y,f(x,y))\vert (x,y)\in D\}\subset M. 
$$

{\bf Definition 5.}
A local group $G$ of transformations on $M$ is called {\it
variational symmetry group for the functional}
$$
{\cal L}[u]=\int \int _{D _0}L(x,y,u^{(2)})dxdy, \leqno(14) 
$$
if for $g_{\varepsilon }\in G,\;g_{\varepsilon }(x,y,u)= (\bar
x,\bar y,\bar u)$, the function $\bar u=\bar f(\bar x,\bar y)= (g\cdot f)(\bar
x,\bar y)$ is defined on $\bar D \subset D _0$
and
$$
\int _{}^{}\int _{\bar D}^{}L(\bar x,\bar y,pr^{(2)} \bar f(\bar x,\bar
y))d\bar x d\bar y= \int \int _{D}L(x,y,pr^{(2)}f(x,y))dxdy. 
$$

The infinitesimal criterion for the variational problem is given by

{\bf Theorem 7.}
{\it A connected group $G$ of transformations acting on 
$M\subset D_0\times U$ is a group of variational symmetries for
the functional (14) if and only if
$$
\hbox{pr}^{(2)}X(L)+L\;Div\xi =0, \leqno(15) 
$$
for $\forall (x,y,u^{(2)})\in M^{(2)}\subset D\times U^{(2)}$ and
for any infinitesimal generator  
$$
X=\zeta (x,y,u){\frac{\partial }{{\partial x}}}+\eta (x,y,u) {\frac{\partial 
}{{\partial y}}}+\phi (x,y,u){\frac{\partial }{{\partial u}}} 
$$
of $G$, where $\xi =(\zeta ,\eta )$ and $Div\xi =D_x\zeta +D_y\eta $.}

{\bf Theorem 8.}
{\it If $G$ is a variational symmetry group of the functional (14),
then $G$ is a symmetry group of Euler-Lagrange equation $E(L)(u)=0$.}

The converse of Theorem 8 is generally false.

{\bf Definition 6.} Let PDE (4) and let $P=(P^1,P^2)$
with $Div\;P=D_xP^1+D_yP^2,$ {\it the total divergence.}
The consequence $Div\;P=0$ of PDE (4) is called {\it conservation law}.
The function $P^1$ is called {\it flow } and $P^2$ is called
{\it conserved density associated to the conservation law.}

By the preceding Definition, there exists a function $Q$ such that 
$$
Div\;P=Q\cdot F.\leqno(16) 
$$

The relation (16) is called {\it the characteristic form of the
conservation law}, and $Q$ is called
{\it the characteristic of the conservation law.}

{\bf Definition 7.}
Let
$$
X=\zeta (x,y,u){\frac{\partial }{{\partial x}}}+\eta (x,y,u){\frac{\partial 
}{{\partial y}}}+\phi (x,y,u){\frac{\partial }{{\partial u}}} 
$$
be a vector field on $M$. The vector field  
$$
X_Q=Q{\frac{\partial }{{\partial u}}},\;\;Q=\phi -\zeta u_x-\eta u_y, 
$$
is called {\it vector field of evolution associated to $X$},
and $Q$ is called {\it the characteristic associated to $X$.}

{\bf Theorem 9 (Noether Theorem).}
{\it Let $G$ be a local Lie group of transformations,
which is a symmetry group of the variational problem (14) and let
$$
X=\zeta (x,y,u){\frac{\partial }{{\partial x}}}+\eta (x,y,u){\frac{\partial 
}{{\partial y}}}+\phi (x,y,u){\frac{\partial }{{\partial u}}} 
$$
the infinitesimal generator of $G$. The characteristic $Q$ of the field $X$
is also a characteristic of the conservation law for the associated 
Euler-Lagrange equation $E(L)(u)=0$.}

There follows the existence of $P=(P^1,P^2)$,
such that
$$
Div\;P=Q\cdot E(L)=0 
$$
to be a conservation law (in the caracteristic form) for the Euler-Lagrange
equation $E(L)=0$.

One proves ([17], 356) that for the Lagrangian $L=L(x,y,u^{(2)})$
we have
$$
P=-(A+L\xi )=-(A^1+L\zeta ,A^2+L\eta )=(P^1,P^2),\;\;A=(A^1,A^2),\leqno(17) 
$$
where  
$$
A^1=Q\cdot E^{(x)}(L)+D_x\left(Q\cdot E^{(xx)}(L)\right)+ {\frac{1}{{2}}}%
D_y\left(Q\cdot E^{(xy)}(L)\right),
\leqno(18)
$$
$$
A^2=Q\cdot E^{(y)}(L)+{\frac{1}{{2}}}D_x\left(Q\cdot E^{(xy)}(L)\right)+
D_y\left(Q\cdot E^{(yy)}(L)\right),
$$
and 
$$
E^{(x)}(L)={\frac{\partial L}{{\partial u_x}}}-2D_x\left({\frac{\partial L}{{%
\partial u_{xx}}}}\right)-D_y\left({\frac{\partial L}{{\partial u_{xy}}}}
\right),
\leqno(19)
$$
$$
E^{(y)}(L)={\frac{\partial L}{{\partial u_y}}}-D_x\left({\frac{\partial L}{{%
\partial u_{xx}}}}\right)-2D_y\left({\frac{\partial L}{{\partial u_{xy}}}}%
\right), 
$$
$$
E^{(xx)}(L)={\frac{\partial L}{{\partial u_{xx}}}},\;\; E^{(xy)}(L)={\frac{%
\partial L}{{\partial u_{xy}}}},\;\; E^{(yy)}(L)={\frac{\partial L}{{%
\partial u_{yy}}}}, 
$$
are {\it Euler operators of superior order.}

\section{Group of variational symmetries of the functional attached
to \c{T}i\c{t}eica PDE. Conservation laws}

We consider the functional  
$$
{\cal L}[u]=\int \int _{D}^{}u\left({\frac{{u_{xx}u_{yy}-u^2_{xy}}}{{%
(xu_x+yu_y-u)^4}}}-\alpha \right)dxdy,\;\;\alpha \in {\bf R}^{*}, \leqno(20) 
$$
where $D$ is a domain in {\bf R$^2$}, $u\in C^{2}(D)$ and the condition
(2) is satisfied for any $(x,y)\in D$.

{\bf Theorem 10.} {\it The Lie algebra of the variational symmetry group of
the functional (20) is described by the vector fields}
$$
Y_1=x{\frac{\partial }{{\partial x}}}-u{\frac{\partial }{{\partial u}}},\;\;
Y_2=y{\frac{\partial }{{\partial y}}}-u{\frac{\partial }{{\partial u}}}, 
\leqno(21) 
$$
$$
Y_3=y{\frac{\partial }{{\partial x}}},\;\;Y_4=x{\frac{\partial }{{\partial y}%
}}. 
$$

{\bf Proof.} According Theorem 8, the vector fields which determine
the Lie algebra of the variational symmetry group are founded
between the vector fields of the Lie algebra of the symmetry group
of the associated Euler-Lagrange equation.
The condition (15) must be verified only for the vector fields in
the Lie algebra (8) of the symmetry group of PDE (1$'$).
One considers 
$$
X=\sum _{i=1}^8C_iX_i,
$$
where $C_i\in {\bf R}$ and $X_i$ are
the infinitesimal generators of the symmetry group $G$
associated to \c{T}i\c{t}eica surfaces PDE.
One determines the real constants $C_i$ such that the relation (15)
is satisfied.
Using the relation (5), the second prolongation of the vector field
$$
X=(C_1x+C_3y+C_4u){\partial \over {\partial x}}+
+(C_5x+C_2y+C_6u){\partial \over {\partial y}}+
(C_7x+C_8y-(C_1+C_2)){\partial \over {\partial u}},
$$
is given by the functions
$$
\Phi ^x=C_7-(2C_1+C_2)u_x-C_5u_y-C_4u^2_x-C_6u_xu_y,
$$
$$
\Phi ^y=C_8-C_3u_x-(C_1+2C_2)u_y-C_4u_xu_y-C_6u^2_y,
$$
$$
\Phi ^{xx}=-(3C_1+C_2)u_{xx}-2C_5u_{xy}-3C_4u_xu_{xx}-C_6u_yu_{xx}-
2C_6u_xu_{xy},
$$
$$
\Phi ^{xy}=-C_3u_{xx}-2(C_1+C_2)u_{xy}-C_5u_{yy}-C_4u_yu_{xx}-2C_6u_yu_{xy}-
$$
$$
-2C_4u_xu_{xy}-C_6u_xu_{yy},
$$
$$
\Phi ^{yy}=-(C_1+3C_2)u_{yy}-2C_3u_{xy}-3C_6u_yu_{yy}-C_4u_xu_{yy}
-2C_4u_yu_{xy}.
$$
Substituting $L$ and $X$ with
$\xi =(C_1x+C_3y+C_4u,C_5x+C_2y+C_6u)$, and
$Div \xi =C_1+C_2+C_4u_x+C_6u_y$ in the relation (15),
after computation, it follows
$$
C_7x+C_8y+C_4uu_x+C_6uu_y=0,
$$
and thus $C_4=C_6=C_7=C_8=0$.
It results that the infinitesimal generator of the
variational symmetry group
for the functional (20) is
$$
X=C_1X_1+C_2X_2+C_3X_3+C_5X_5.
$$
Denote $Y_1=X_1,\;Y_2=X_2,\;Y_3=X_3$ and $Y_4=X_5$.

{\bf Proposition 2.}
{\it For the vector field
$$
-Y_3=-y{\frac{\partial }{{\partial x}}}, 
$$
the flow and respectively the
conserved density of the conservation law are} 
$$
P^1=-\alpha yu+{\frac{u_x}{{(xu_x+yu_y-u)^4}}}(u_{xy}(yu_y-u)-yu_xu_{yy}),
\leqno(22) 
$$
$$
P^2=-{\frac{u_x}{{(xu_x+yu_y-u)^4}}}(u_{xx}(yu_y-u)-yu_xu_{xy}). 
$$

{\bf Proof.}
The caracteristic associated to the vector field $-Y_3$ is
$$
Q=yu_y. 
$$
Replacing in the relations (18), we obtain 
$$
A^1={\frac{{uy(u^2_{xy}-u_{xx}u_{yy})}}{{(xu_x+yu_y-u)^4}}}+ {\frac{{%
u_{xy}(uu_x-yu_y)}}{{(xu_x+yu_y-u)^4}}}+ {\frac{{yu^2_xu_{yy}}}{{%
(xu_x+yu_y-u)^4}}}, 
$$
$$
A^2={\frac{u_x}{{(xu_x+yu_y-u)^4}}}(u_{xx}(yu_y-u)-yu_xu_{xy}). 
$$

Introducing $\xi =(-y,0)$ in the relations (17) it follows that the
functions $P^1,P^2$ have the form (22).

Analogously one determines the conservation laws corresponding to the
characteristics of the vector fields (21).

\section{Strong/weak symmetry group}

The symmetry group introduced in the Definition 2 is called {\it strong
symmetry group}.

{\bf Definition 8.}
{\it The weak symmetry group of PDE} (4)
is a group of transformations acting on $M\subset D\times U$
and which satisfies only the condition (b) in the Definition 2 of the
strong symmetry group.

Consequently a weak symmetry group did not transforms solutions of PDE
into its solutions.

{\bf Proposition 3.}
{\it Let $G$ be a connected Lie group of transformations on $M$,
with infinitesimal generators $X_1,...,X_s$. Let $Q^1,...,Q^s$
be the characteristics associated to these vector fields. Then any 
$G$-invariant function $u=f(x,y)$ must satisfy the system of equations}
$$
Q^k(x,y,u^{(1)})=0,\;\;k=1,...,s. \leqno(23) 
$$

Any $G$-invariant solution $u=f(x,y)$ of PDE (4) is also a solution of the
system (23), and hence of the system 
$$
\left\{
\begin{array}{ccl}
F(x,y,u^{(2)}) & = & 0 \\ 
Q^k(x,y,u^{(1)}) & = & 0, \;\;k=1,...,s.
\end{array}
\right. \leqno(24) 
$$
The converse is true only for the case in which $G$ is a strong symmetry
group.

One proves ([19])

{\bf Theorem 11.}
{\it Let $G$ be a group of transformations acting on
$M\subset D\times U$ and (4) a PDE of order two defined on $D$.
Then $G$ is always a symmetry group of the system (24) and hence always a
weak symmetry group.}

Every $s$-parameter subgroup $H$ of the strong symmetry group $G$ (8)
determines a family of group-invariant solutions. The problem of
classification of the group-invariant solutions is reduced to the
problem of classification of Lie subalgebras of the Lie algebra {\bf g}
of the group $G$ ([14], 186).
For the 1-dimensional subalgebras one considers a general element $X$
and this can be simplified as much as possible, using the adjoint
transformations.

We shall determine some solutions of PDE (1$'$) which are invariant with
respect to the strong symmetry group $G$.

{\bf Remarks.}

1) The finding of the adjoint representation $Ad\;G$ of the  Lie group $G$,
can be realised using the Lie series
$$
Ad(exp(\varepsilon X)Y)=\sum \limits _{n=0}^{\infty }{\frac{{\varepsilon ^n} 
}{{n!}}}(adX)^n(Y)=Y-\varepsilon [X,Y]+{\frac{\varepsilon ^2}{{2}}}
[X,[X,Y]]-... \leqno(25) 
$$

2) If $u=f(x,y)$ is a solution of PDE (1$'$), then the following functions
$$
u^{(1)}=e^{-\varepsilon }f(xe^{-\varepsilon },y),\;\;
u^{(2)}=e^{-\varepsilon }f(x,ye^{-\varepsilon }),\;\;%
u^{(3)}=f(x-\varepsilon y,y),
$$
$$
u^{(4)}=f(x-\varepsilon u^{(4)},y),\;\;
u^{(5)}=f(x, y-\varepsilon x),\;\; u^{(6)}=f(x, y-\varepsilon u^{(6)}),\;\;%
$$
$$
u^{(7)}=f(x,y)+\varepsilon x,\;\; u^{(8)}=f(x,y)+\varepsilon
y,\;\;\varepsilon \in {\bf R},
$$
are also solutions of the equation since every 1-parameter subgroup
$G_i$ generated by $X_i,\;i=1,...,8$, is a symmetry group.

3) The adjoint representation $Ad\;G$ of the Lie group $G$
which invariates the \c{T}i\c{t}eica equation, is determined using
the Lie series (25). This way we obtain

\vspace{0,5cm}

\centerline{\begin{tabular}{|c|c|c|c|c|}
\hline
$\scriptstyle{Ad}$ & $\scriptstyle{X_1}$ & $\scriptstyle{X_2}$ &
$\scriptstyle{X_3}$ & $\scriptstyle{X_4}$ \\
\hline
$\scriptstyle{X_1}$ & $\scriptstyle{X_1}$ & $\scriptstyle{X_2}$ &
$\scriptstyle{e^{\varepsilon }X_3}$ & $\scriptstyle{e^{2\varepsilon }X_4}$ \\
\hline
$\scriptstyle{X_2}$ & $\scriptstyle{X_1}$ & $\scriptstyle{X_2}$ &
$\scriptstyle{e^{-\varepsilon }X_3}$ & $\scriptstyle{e^{\varepsilon }X_4}$ \\
\hline
$\scriptstyle{X_3}$ & $\scriptstyle{X_1-\varepsilon X_3}$ &
$\scriptstyle{X_2+\varepsilon X_3}$ & $\scriptstyle{X_3}$ &
$\scriptstyle{X_4}$ \\
\hline
$\scriptstyle{X_4}$ & $\scriptstyle{X_1-2\varepsilon X_4}$ &
$\scriptstyle{X_2-\varepsilon X_4}$ & $\scriptstyle{X_3}$ &
$\scriptstyle{X_4}$ \\
\hline
$\scriptstyle{X_5}$ & $\scriptstyle{X_1+\varepsilon X_5}$ &
$\scriptstyle{X_2-2\varepsilon X_5}$ &
$\scriptstyle{X_3-\varepsilon (X_1-X_2)-{\varepsilon }^2X_5}$ &
$\scriptstyle{X_4+\varepsilon X_6}$ \\
\hline
$\scriptstyle{X_6}$ & $\scriptstyle{X_1-\varepsilon X_6}$ &
$\scriptstyle{X_2-2{\varepsilon }X_6}$ &
$\scriptstyle{X_3-\varepsilon X_4}$ & $\scriptstyle{X_4}$ \\
\hline
$\scriptstyle{X_7}$ & $\scriptstyle{X_1+2\varepsilon X_7}$ &
$\scriptstyle{X_2+\varepsilon X_7}$ &
$\scriptstyle{X_3+\varepsilon X_8}$ & $\scriptstyle{X_4-{\varepsilon }X_1-
{\varepsilon }^2X_7}$  \\
\hline
$\scriptstyle{X_8}$ & $\scriptstyle{X_1+\varepsilon X_8}$ &
$\scriptstyle{X_2+2\varepsilon X_8}$ & $\scriptstyle{X_3}$ &
$\scriptstyle{X_4-\varepsilon X_3}$ \\
\hline
\end{tabular}}

\par

\vspace{0,5cm}
\centerline{\begin{tabular}{|c|c|c|c|c|}
\hline
$\scriptstyle{Ad}$ & $\scriptstyle{X_5}$ & $\scriptstyle{X_6}$ &
$\scriptstyle{X_7}$ & $\scriptstyle{X_8}$ \\
\hline
$\scriptstyle{X_1}$ & $\scriptstyle{e^{-\varepsilon }X_5}$ &
$\scriptstyle{e^{\varepsilon }X_6}$ & $\scriptstyle{e^{-2\varepsilon }X_7}$ &
$\scriptstyle{e^{-\varepsilon }X_8}$ \\
\hline
$\scriptstyle{X_2}$ & $\scriptstyle{e^{\varepsilon }X_5}$ &
$\scriptstyle{e^{2\varepsilon }X_6}$ & $\scriptstyle{e^{-\varepsilon }X_7}$ &
$\scriptstyle{e^{-2\varepsilon }X_8}$ \\
\hline
$\scriptstyle{X_3}$ & $\scriptstyle{X_5-\varepsilon (X_1-X_2)-{\varepsilon }
^2X_3}$ & $\scriptstyle{X_6+\varepsilon X_4}$ & $\scriptstyle{X_7-\varepsilon
X_3}$ & $\scriptstyle{X_8}$ \\
\hline
$\scriptstyle{X_4}$ & $\scriptstyle{X_5-\varepsilon X_6}$ &
$\scriptstyle{X_6}$ & $\scriptstyle{X_7+\varepsilon X_1-
{\varepsilon }^2X_4}$ & $\scriptstyle{X_8+\varepsilon X_3}$ \\
\hline
$\scriptstyle{X_5}$ & $\scriptstyle{X_5}$ & $\scriptstyle{X_6}$ &
$\scriptstyle{X_7}$ & $\scriptstyle{X_8-\varepsilon X_7}$ \\
\hline
$\scriptstyle{X_6}$ & $\scriptstyle{X_5}$ & $\scriptstyle{X_6}$ &
$\scriptstyle{X_7+\varepsilon X_5}$ & $\scriptstyle{X_8+\varepsilon X_2-
{\varepsilon }^2X_6}$ \\
\hline
$\scriptstyle{X_7}$ & $\scriptstyle{X_5}$ & $\scriptstyle{X_6-\varepsilon X_5}$
& $\scriptstyle{X_7}$ & $\scriptstyle{X_8}$ \\
\hline
$\scriptstyle{X_8}$ & $\scriptstyle{X_5+\varepsilon X_7}$ &
$\scriptstyle{X_6-\varepsilon X_2-{\varepsilon }^2X_8}$ &$\scriptstyle{X_7}$
& $\scriptstyle{X_8}$ \\
\hline
\end{tabular}}

\vspace{0,5cm}

Finally, we determine some group-invariant solutions of the equation
(1$'$), corresponding to 1-dimensional subalgebras generated by
$X_1-X_2,\;\;X_5-X_3$.

a) For the vector field
$$
X_1-X_2=x{\partial \over {\partial x}}-y{\partial \over {\partial y}},
$$
one looks for solutions of the form $u=\varphi (xy)$.
In this case the PDE (1$'$) becomes
$$
2t\varphi '\varphi ''+\varphi '^2+\alpha (2t\varphi '-\varphi )^4=0,
$$
where $t=xy.$
This DE admits particular solutions of the form $\phi (t)=t^p$,
with the condition imposed in (2).
For $p=-1$ and $\alpha ={1\over {27}}$, we obtain 
$$
u(x,y)={1\over {xy}}
$$
as a solution of PDE (1$'$). According to the preceding remark 2, it follows that
$$
u(x,y)={1\over {xy}}+\varepsilon x,\;\;u(x,y)={1\over {xy}}+
\varepsilon y,\;\;
u(x,y)={1\over {(x-\varepsilon _1y)y}}+\varepsilon _2x,\;\;\varepsilon ,
\varepsilon _1,\varepsilon _2\in {\bf R},
$$
are also solutions.

Other particular solution of the preceding DE is
$\varphi (t)=\sqrt {1+at}$, for
$a^2+4\alpha $=0. For $\alpha <0$, it follows the solution 
$$
u(x,y)=\sqrt{1+axy},\;\;a\in {\bf R}^{*},
$$
of PDE (1$'$). According to remark 2, the functions
$$
u(x,y)=\sqrt{1+axy}+\varepsilon x,\;\;
u(x,y)=\sqrt{1+a(x-\varepsilon _1y)y}+\varepsilon _2x,\;\;\varepsilon ,
\varepsilon _1,\varepsilon _2\in {\bf R},
$$
are also solutions of PDE (1$'$).

b) For the vector field 
$$
X_5-X_3=-y{\partial \over {\partial x}}+x{\partial \over {\partial y}},
$$
one looks for solutions of the form $u=\varphi (r),$ where $r=\sqrt{x^2+y^2}.$
Replacing in the PDE (1$'$) we obtain the DE
$$
{1\over r}\varphi '\varphi ''=\alpha (r\varphi '-\varphi )^4.
$$
This DE admits particular solutions of the form $\phi (r)=r^p$.
For $p=-2$ and $\alpha =-{4\over {27}}$, it follows 
$$
u(x,y)={1\over {x^2+y^2}}
$$
as solution of PDE (1$'$). According to remark 2, the functions
$$
u(x,y)={1\over {x^2+y^2}}+\varepsilon x,\;\;
u(x,y)={1\over {{(x-\varepsilon _1y)}^2+y^2}}+\varepsilon _2x,\;\;
\varepsilon ,\varepsilon _1,\varepsilon _2\in {\bf R},
$$
are also solutions of PDE (1$'$).
The DE admits also a particular solution of the form $\varphi (r)=\sqrt{1+ar^2},
\;a\in {\bf R}^{*}$ for $\alpha =a^2$.
If $\alpha >0$, then it follows the implicit solution
$$
u^2+a(x^2+y^2)=1,\;\;u>0,
$$
of PDE (1$'$), and according to remark 2, the equations
$$
u^2+a({(x-\varepsilon y)}^2+y^2)=1,\;\;u>0,
$$
$${(u-\varepsilon x)}^2+a(x^2+y^2)=1,\;\;u-\varepsilon x>0,\;\;
\varepsilon \in {\bf R},
$$
define also solutions of PDE (1$'$).

Now we refer to weak symmetry groups and the corresponding
solutions of PDE (1$'$).

a) Let
$$
X=-x^2y{\partial \over {\partial x}}+{\partial \over {\partial u}}.
$$
We obtain $C_1=u-{1\over {xy}},\;\;C_2=y$. Hence
$$
u(x,y)={1\over {xy}}+h(y).
$$
Replacing in the PDE (1$'$), it follows the DE
$$
{3\over {x^4y^4}}+{2h''\over {x^3y}}=\alpha {\left(-{3\over {xy}}+
yh'-h\right)}^4.
$$
As $h=h(y)$, by identification we deduce
$h''=0,\;\;yh'-h=0,\;\;\alpha ={1\over {27}}$,
hence $h(y)=Cy,\;C\in {\bf R}$. Consequently
$$
u(x,y)={1\over {xy}}+Cy,\;\;C\in {\bf R}
$$
is a solution of PDE (1$'$).

b) Let
$$
X=2ux{\partial \over {\partial x}}+(u^2-1){\partial \over {\partial u}}.
$$
Since $C_1={{u^2-1}\over {x}},\;\;C_2=y$, it follows 
$$
u(x,y)=\sqrt{1+xh(y)},
$$
for $u>0$.

Replacing in PDE (1$'$) we obtain $h(y)=Cy,\;C\in ${\bf R}$^{*}$.
The corresponding solution of PDE (1$'$) is
$$
u(x,y)=\sqrt{1+Cxy},\;\;u>0,\;\;C\in {\bf R}^{*},\;\;1+Cxy>0.
$$

\par
\centerline{University "Politehnica" of Bucharest}
\centerline{Departament of Mathematics I}
\centerline{Splaiul Independen\c{t}ei 313}
\centerline{77206 Bucharest Romania}
\centerline{e-mail:udriste@mathem.pub.ro}
\centerline{e-mail:nbila@mathem.pub.ro}

\end{document}